\documentclass{amsart}

\usepackage{amsfonts,amssymb,amsmath,amsthm,url}

\def\tfr#1/#2{{#1}/{#2}}				
\def\dfr#1/#2{\displaystyle\frac{#1}{#2}}	
\def\fr#1/#2{\frac{#1}{#2}}				

\def\e#1{e^{#1}}						
\def\Cplx{\mathbb{C}}

\usepackage[pdftex]{graphicx}

\def\textfig#1#2#3#4#5{
\noindent
\begin{minipage}[h]{0.5\textwidth}
	{#1}
\end{minipage}
\hfill
\begin{minipage}[h]{0.5\textwidth}
	\scalebox{#2}{%
\begingroup
  \setlength{\unitlength}{1bp}
  \begin{picture}(#3,#4)
    \put(0,0){\includegraphics{#5}}
  \end{picture}
\endgroup
}
\end{minipage}
}

\author{Jacques G\'elinas}

\address{Ottawa, Canada}

\email{jacquesg00@hotmail.com}

\thanks{This work was done in 2018 while the author was a retired mathematician}

\begin{document}




\title[Bernoulli numbers]{Definitions, notations and proofs for Bernoulli numbers}

\begin{abstract}

This is a collection of definitions, notations and proofs for the Bernoulli numbers $B_n$ appearing in formulas for the sum of integer powers, some of which can be found scattered in the large related historical literature in French, English and German.  We provide elementary proofs for the original convention with ${\mathcal B}_1=1/2$ and also for the current convention with $B_1=-1/2$, using only the binomial theorem and the concise Blissard symbolic notation.

\end{abstract}

\maketitle


\section{Historical definition and notations}
The Bernoulli numbers are an infinite sequence of rational coefficients appearing in the polynomials expressing the sum of the powers of the first $n$ natural numbers as a function of $n$. Jacob Bernoulli (1655--1705) was the first to define them in his book ``Ars Conjectandi" published posthumously in 1713.
This book, one of the first on probability, was translated with explanations into English in 1795, into German in 1899, and again into English in 2006. Its second part deals with counting permutations and combinations, and the sums of powers are included there on pages 96--97.

\textfig{
Bernoulli states that in a table of combinations, figurate numbers, or binomial coefficients, each element is equal to the sum of the numbers above it in the previous column~:
$$
	\sum_{k=1}^{n} \binom{k-1}{p} = \binom{n}{p+1},\qquad(n\ge1,p\ge0).
 $$
This is indeed easily proven by Bernoulli's own telescoping summation method, or by counting combinations, and gives the convenient recurrence relation explained below for the sums of powers of integers which he denoted by $\int n^p,\,p\ge0$.
}{1.1}{200}{150}{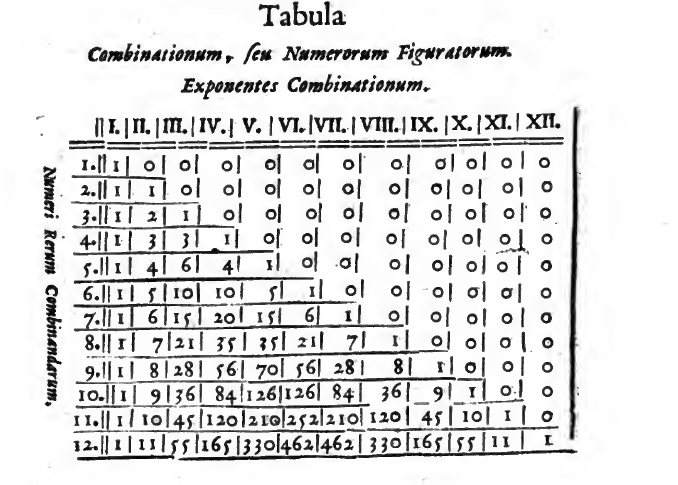}
\smallskip
\textfig{If $p=1$ for example, then $\binom{k-1}{1}=k-1$ and
\begin{align*}
	&\sum_{k=1}^{n} \binom{k-1}{1} = \int n^1 - \int n^0 = \binom{n}{2}
\\
	&\qquad\qquad	\implies \int n^1 = \fr {n^2}/2 + \fr n/2;
\end{align*}
if $p=2$, then $\binom{k-1}{2}=\tfr k^2/2- \tfr 3k/2 + 1$ and
\begin{align*}
	&\sum_{k=1}^{n} \binom{k-1}{2} = \fr 1/2\int n^2 -\fr 3/2\int n^1 + \int n^0 = \binom{n}{3}
\\
	&\qquad\qquad	\implies \int n^2 = \fr {n^3}/3 + \fr {n^2}/2 + \fr {n}/6.
\end{align*}
The derivation of the sum of the cubes is also explained (up to $p=10$ in the 1795 translation).
}{1.1}{200}{150}{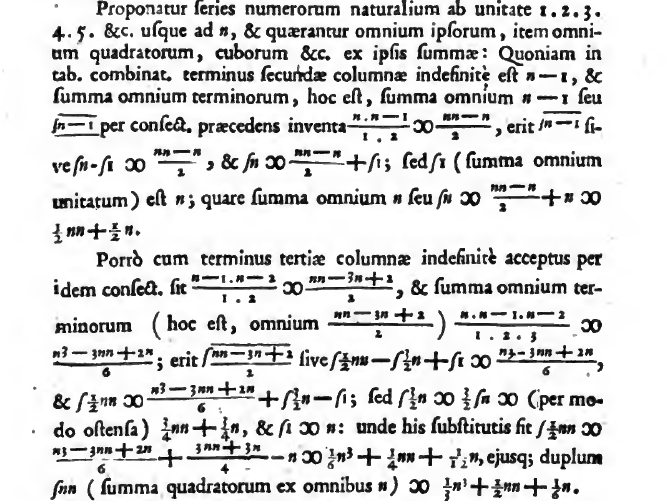}

\smallskip
\textfig{
Bernoulli boasted that he could thus derive the sums of the first six powers on a single page, a result obtained previously with difficulty by Ismael Bulliadus (1631--1694). In fact he provided a table of the first ten formulas, so that the exact sum of the tenth powers of the the first $1000$ integers could be computed in ``half the quarter of an hour". Much more importantly, he gave an explicit formula for $\int n^p$, featuring the numbers $\tfr 1/6, -\tfr 1/{30}, \tfr 1/{42}, -\tfr 1/{30}, \tfr 5/{66},\ldots$ as coefficients of $n$ in the successive $\int n^{2p}, p\ge1$. Finally, he indicated that other numbers in the list could be computed successively by increasing $p$ and using $n=1$ (since then
 $\int n^{2p}=1$). However, he provided no formal proof of this rule or of the formula itself.
}{1.1}{200}{150}{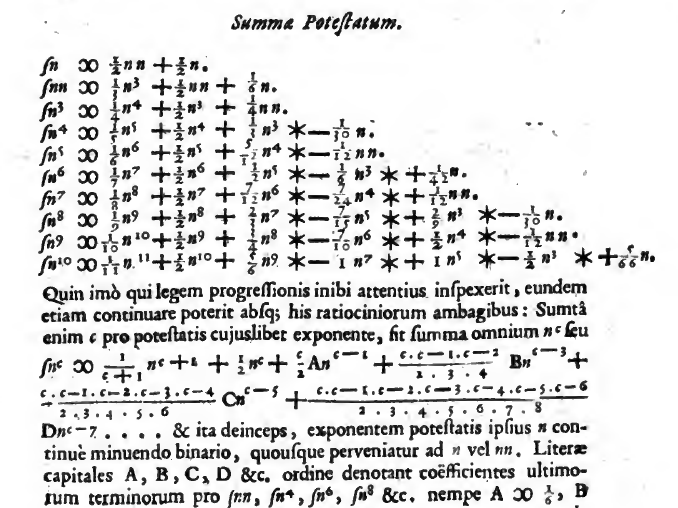}

In 1730, Abraham De Moivre expressed the rule of Jacob Bernoulli as a recurrence relation for the computation of what he called ``Bernoulli numbers",
$B_2=\tfr 1/6$, $B_4=-\tfr 1/{30}$ and so forth~:
$$
	\binom{2m+1}{1}B_{2m}+\binom{2m+1}{3}B_{2m-2}+\ldots+\binom{2m+1}{2m-1}B_2 = \fr 2m-1/2.
 $$
This was included in the 1800 treatise on finite differences of Lacroix, with different signs and indices (see below right). De Moivre had used it to reformulate and prove an asymptotic relation communicated to him by Stirling and featuring a series shown later by Bayes to be divergent \cite{Gelinas:Stirling}~:
$$
\log (n-1)! \approx \left(n-\fr 1/2\right)\log n - n + \log\sqrt{2\pi}
	+ \sum_{k=1}^\infty \fr {B_{2k}}/{2k(2k-1)n^{2k-1}}.
 $$

\textfig{
Next, the publication in 1755 of ``Institutiones Calculi Differentialis" by Euler was a revolution in the early history of these Bernoulli numbers. He proved that these numbers are the Maclaurin coefficients of the even ``generating function" $(\tfr x/2) \coth(\tfr x/2)$, that their signs alternate, and that they appear in the Euler-Maclaurin summation formula and in the values of the infinite series $\zeta(2p)=\sum\tfr 1/{n^{2p}}$.
He also computed the values of the next 10 numbers $B_{12},\ldots,B_{30}$ (\S V.132), via a binomial convolution formula obtained from the differential equation of $y=\cot x$($y'+y^2+1=0$) whose series coefficients after $\tfr 1/x$ are all negative and prove the alternating sign rule (\S V.119,\cite{Gessel:2005}).
}{0.4}{475}{420}{J_Gelinas-2019-Bernoulli_numbers-Lacroix}

The general Bernoulli formula can be displayed as follows, using the convention adopted early in the 19th century of summing only the first $n-1$ powers while starting at $0$ so that $S_0(n)=n$ when we agree to define $0^0:=1$. Alternating signs and missing terms stand out clearly.
\begin{align*}
	S_p(n) &:= 0^p+1^p+2^p+\ldots+(n-1)^p			
\\
	&= \fr 1/{p+1}\left[ \binom{p+1}{0} n^{p+1} - \fr 1/2 \binom{p+1}{1}n^p
		+ \fr 1/6 \binom{p+1}{2}n^{p-1} - \fr 1/{30} \binom{p+1}{4}n^{p-3} + \ldots
	\right].
\end{align*}
The current \emph{even index} notation, already used by De Morgan in 1836, includes all the successive coefficients of the polynomials, $B_0=1$, $B_1=-\tfr 1/2$, $B_2=\tfr 1/6$, $B_3=0$, $B_4=-\tfr 1/{30}$, $B_5=0$, and so forth, allowing the use of the summation symbol as in
\begin{equation} \label{Bernsum}
	S_p(n) = \fr 1/{p+1}\,\sum_{k=0}^p \binom{p+1}{k} B_k n^{p+1-k},	\qquad(n\ge0,p\ge0).
\end{equation}

	We can simplify this further with the \emph{representative} symbolic notation introduced by the Rev. John Blissard in a series of articles on ``Generic Equations" published during the 1860s in the Quarterly Journal of Mathematics edited by Arthur Cayley and J.J. Sylvester. This notation was quickly adopted by Lucas (1870s) and Ces\`aro (1880s) and has become standard (N\"orlund 1924). 

First, any real polynomial of degree $n$ can obviously be written with binomial factors as
$$
	P(x) = \sum_{k=0}^n \binom{n}{k} a_k x^{n-k}.
 $$

Secondly, if we introduce a dummy variable $A$ as a \emph{representative} of the coefficients $a_k$ and agree to \emph{downgrade} the exponents of $A$ to indices \emph{after expansion} via the binomial theorem (including $A^0\to a_0$), then we can write symbolically
$$
	P(x) = (A + x)^n.
  $$
It is also easily verified that $n(A+x)^{n-1}$ is a symbolic representation of $P'(x)$, as it should.

	If B is the representative of the Bernoulli numbers, (\ref{Bernsum}) thus becomes
\begin{equation} \label{Bernbinom}
	S_p(n) = \fr {(B+n)^{p+1} - B_{p+1}}/{p+1} \qquad(n\ge0,p\ge0),
\end{equation}
or, equivalently,
\begin{equation} \label{Bernint}
	S_p(n) = \int_0^n (B+x)^p\,dx, \qquad(n\ge0,p\ge0).
\end{equation}

	If we accept the validity of (\ref{Bernbinom}), we also obtain immediately by differentation
\begin{equation} \label{Bernderiv}
	S'_p(n) = (B+n)^p = p S_{p-1}(n) + B_p, \qquad(n\ge0,p\ge1).
\end{equation}
This suggests a fast method \cite{Prouhet:1864} for computing the numeric coefficients of the polynomials $S_p(n)$, increasing successively the exponent $p$ by unity from $1$~: add to the anti-derivative of $pS_{p-1}(n)$ a last linear term $B_p n$, choosing $B_p$ so that $S_p(1)=0$. Thus, from $S_0(n)=n=B_0n$,
$$
\begin{array}{ll}
	S_1(n) =\int n\,dn + B_1n = \fr n^2/2 + B_1n
			&\implies B_1=-\fr 1/2,
\\[0.5em]
	S_2(n) =\int (n^2-n)\,dn + B_2n  = \fr n^3/3 - \fr n^2/2 + B_2n
			 &\implies B_2 = \fr 1/2 - \fr1/3 = \fr 1/6,
\\[0.5em]
	S_3(n) =\int \left(n^3-\fr 3/2n^2+\fr n/2\right)\,dn + B_3n
			= \fr n^4/4 - \fr n^3/2 + \fr n^2/4 + B_3n &\implies B_3=-\fr 1/4 + \fr 1/2 - \fr 1/4 = 0.
\end{array}
 $$
This also works for the original Bernoulli $\int n^p$, but the sum of the coefficients must now be $1$ instead of $0$, since adding $n^p$ to $S_p(n)$ changes the sign of the  second term $-\tfr n^p/2$ on the right of (\ref{Bernsum}).

\def\B^+{\mathcal B}

The usual convention of summing up to $n-1$ can be justified by an asymmetric property of finite differences and the convenience of having all the Bernoulli polynomials $B_p(x):=(B+x)^p$ symmetric or antisymmetric about $x=\tfr 1/2$. On the other hand, as noted by Lucas and Ces\`aro, symbolic formulas are simpler and more natural with $\B^+_1=+\tfr 1/2$, such as $\gamma=-\log \B^+$ (Ces\`aro, 1880) or $(s-1)\zeta(s) = \B^+^{1-s}, \, s\in \Cplx$ (Yiping Yu, 2012). Adopting the usual notation with $B_1=-\tfr 1/2$ often implies the addition of an extra sign factor only needed for $p=1$ as in $(-1)^pB_p=-p\zeta(1-p),\,p\ge0$.

	Throughout the 19th century, only the positive coefficients with consecutive indices in the Bernoulli formula were most often considered, and usually defined via the Euler identity
$$
	\fr 1/{\e{x}-1} = \fr 1/x - \fr 1/2 + \sum_{n=1}^\infty (-1)^{n-1} \fr {B^{\text{old}}_n}/{(2n)!} x^{2n-1}.
 $$
With this 19th century notation, the Bernoulli rule (De Moivre recurrence) was written as
$$
	\binom{2p+1}{1}B_p^{\text{old}} - \binom{2p+1}{3}B_{p-1}^{\text{old}} + 
	\binom{2p+1}{5}B_{p-2}^{\text{old}} - \ldots	+ (-1)^p\fr {2p-1}/2 = 0,
	\qquad(p=0,1,2,\ldots).
 $$
The Blissard representative notation gives us more concise and mnemonic recurrence formulas~:
\begin{equation} \label{Bernrecur}
	(B+1)^{p+1}-B_{p+1}=0^p \iff (B+1)^p=(-B)^p,	\qquad(p=0,1,2,\ldots).
\end{equation}


\section{Proofs that Jacob Bernoulli could have provided}
	In his 1893 treatise ``Vorlesungen {\"u}ber die Bernoullischen Zahlen", Louis Saaltsch{\"u}tz noted that he was not aware, two centuries later, of anyone using the original definition of Bernoulli to prove the sum of powers formula, given the availability of the clever and powerful generating function approach of Euler. Saaltsch{\"u}tz pointed out the following particular properties as needing justification:
\begin{enumerate}
\item{} The $\int n^p$ can be expressed by polynomials of degree $p+1$ in the variable $n$.
\item{} The $B_{2k}$ are constants invariant with respect to the exponent $p$.
\item{} The $B_{2k+1}$ are zero for $k>0$ and all terms in $\int n^p - \tfr {n^p}/2$ have the parity of $p+1$
if $p>0$.
\item{} The $B_{2k}$ alternate in sign for $k>0$.
\end{enumerate}

	This section fills the gap in a shorter fashion than Saaltsch{\"u}tz did in the first chapter of his book.
Let us start by assembling the following table, using current notations with $B_1=-\tfr 1/2$.
$$
\begin{array}{|l|l|l|} \hline
		& \text{Property}	& \text{Equation}
\\[0.5em] \hline
	\text{(A)} & \text{Definition}	&	S_p(n) := \sum_{k=0}^{n-1} k^p,\qquad(n\ge1,p\ge0,0^0:=1)	\phantom{\bigg[}
\\[0.5em]
	\text{(B)}	& S\text{--$p$--recurrence} &	\sum_{k=0}^p\binom{p+1}{k}S_k(n) = n^{p+1}
\\[0.5em]
	\text{(C)} & \text{First $S_p$ polynomials} & 	S_0(n) = n; S_1(n) =\tfr{n(n-1)}/2
\\[0.5em]
	\text{(D)}	& \text{Degree}		&	S_p(n) \text{ is a polynomial of degree } p+1
\\[0.5em]
	\text{(E)}	& S\text{--$n$--recurrence} &	S_p(n+1) = S_p(n) + n^p
\\[0.5em]
	\text{(F)}	& \text{First $S_p$ values} &	S_p(1)=0^p; S_p(0)=0; S_p(-1)=(-1)^{p+1}
\\[0.5em]
	\text{(G)}	& \text{Notation}		&	S_p(n) =\ldots+B_p n, \text{ or } B_p:=\tfr S_p(n)/n\Big|_{n=0}
\\[0.5em]
	\text{(H)}	& B\text{--recurrence}	&	(B+1)^{p+1}-B_{p+1} := \sum_{k=0}^p\binom{p+1}{k}B_k = 0^p
\\[0.5em]
	\text{(I)}	& \text{First numbers}	&	B_0=1; B_1=-\fr 1/2; B_2=\fr 1/6; B_3=0; B_4=-\fr 1/{30}
\\[0.5em]
	\text{(J)}	& \text{Sums of powers}	&	S_p(n) = \dfr {(B+n)^{p+1}-B_{p+1}}/{p+1}
									:= \dfr {\sum_{k=0}^p\binom{p+1}{k}B_kn^{p+1-k}}/{p+1}
\\[0.5em]
	\text{(K)}	& \text{Missing terms}	&	B_{2k+1}=0, \text{ for } k>0
\\[0.5em]
	\text{(L)}	& \text{Parity}			&	S_p(n) + \tfr n^p/2 \text{ has the parity of } p+1 \text{ if } p>0	\phantom{\dfr 1/2}
\\
\hline
\end{array}
 $$

	The $S$--$p$--recurrence (B), discovered by Pascal around 1654, comes from a telescoping sum~:
\begin{align*}
	n^{p+1}-0^{p+1}&= \sum_{k=0}^{n-1} \left[(k+1)^{p+1}-k^{p+1}\right]
				= \sum_{k=0}^{n-1}\sum_{j=0}^p\binom{p+1}{j}k^j
\\
			&= \sum_{j=0}^p\binom{p+1}{j}\sum_{k=0}^{n-1} k^j
			 = \sum_{j=0}^p\binom{p+1}{j}S_j(n).
\end{align*}

	By induction using (B) and (C), $S_p(n)$ is therefore a polynomial of degree $p+1$ without constant term in the variable $n$, proving (D) and also giving $S_p(0)=0$.

	The $S$--$n$--recurrence (E) is equivalent to the definition (A) and remains valid for $n=-1$, since the difference of the two sides of its equation is a polynomial vanishing for all $n>0$. This yields $S_p(-1)=(-1)^{p+1}$ from $S_p(0)=0$; and we have $S_p(1)=1$ by the definition (A), proving (F).

	Since $S_p(0)= 0$, $B_p$ is well defined in (G). Dividing the $S$--$p$--recurrence (B) by $n$ and then setting $n=0$ yields the $B$--recurrence (H), since $B_p$ is the coefficient of $n$ in $S_p(n)$ by (G). The $B$--recurrence (H) has a unique solution, which shows the invariance of the $B_k$ with respect to the exponent $p$.
	
	The values in (I) are obtained from (C) or from (H).

\newpage

	We will now deduce (J) from (H), by induction on $n$ for a fixed $p\ge0$. 
First, (J) is verified for all $p\ge0$ when $n=0$ since both sides vanish by (F). Next, if (J) is true for a fixed $p\ge0$ and a certain $n\ge0$, then (J) is also true for $n+1$, from the $B$--recurrence (H) and the $S$--$n$--recurrence (E)~:
\begin{align*}
  (B+n+1)^{p+1} - B_{p+1} &= (B+n)^{p+1} - B_{p+1} 
					+ \sum_{k=0}^{p+1}\binom{p+1}{k}\left[(B+1)^k-B_k\right]n^{p+1-k}
\\
	&= (p+1)S_p(n) + \sum_{k=0}^{p+1}\binom{p+1}{k}0^{k-1} n^{p+1-k}
\\
	&= (p+1)(S_p(n)+n^p) = (p+1)S_p(n+1).
\end{align*}

Next, repeated applications of the $S$-$n$-recurrence (E) starting with negative $n$, given the value of $S_p(0)=0$ from (F), prove the parity property (L), from which (K) follows by (G).
\begin{align*}
	S_p(-n) + \fr (-n)^p/2 &= S_p(-n+1) -(-n)^p +  \fr (-n)^p/2 = 0  - (-1)^p - \ldots -\fr(-n)^p/2
\\
	&= (-1)^{p+1}\left(0^p+1^p+\ldots+\fr n^p/2\right) = (-1)^{p+1}\left(S_p(n)+\fr n^p/2\right).
\end{align*}

	This concludes the justification of our table, and we can boast of having proven an infinity of formulas for the sums of powers in one page and a quarter, using only the binomial theorem.

The alternating sign property can be proven by the following Abel integral formula, equivalent to the Euler equation $(-1)^{n-1}B_{2n}/{(2n)!}=\tfr {2\zeta(2n)}/{(2\pi)^{2n}}$,
$$
	\cos Bx :=\sum_{k=0}^\infty (-1)^n B_{2n} \fr {x^{2n}}/{(2n)!}
		= \fr x/2\cot \fr x/2 = 1 - 2 x \int_0^\infty \fr {\sinh xt}/{\e{2\pi t}-1}\,dt,
	\qquad(|x|<2\pi).
 $$

\def\ml{\mathcal}
\def\B^+{\ml B}
\def\S^+{\ml S}
\def\ml{\mathfrak}
\def\B^+{\ml B}
\def\S^+{\mathcal S^+}
\def\S^+{\mathcal T}
\def\S^+{T}

	Next, we convert the previous table to the conventions used by Jacob Bernoulli, Euler, and Saaltsch{\"u}tz, but we will replace $\int n^p$ by $\S^+_p(n)$ and $B_p$ by $\B^+_p$, including again null coefficients. 
$$
\begin{array}{|l|l|l|} \hline
		& \text{Property}	& \text{Equation}
\\[0.5em] \hline
	(\alpha)	& \text{Definition}	&	\S^+_p(n) := \sum_{k=1}^n k^p,\qquad(n\ge1,p\ge0)	\phantom{\bigg[}
\\[0.5em]
	(\beta)	& \S^+\text{--$p$--recurrence} &	\binom{\S^+-1}{p} = \binom{n}{p+1}
\\[0.5em]
	(\gamma)	& \text{First $\S^+$ polynomials} & 	\S^+_0(n) = n; \S^+_1(n) =\tfr{n(n+1)}/2
\\[0.5em]
	(\delta)	& \text{Degree}		&	\S^+_p(n) \text{ is a polynomial of degree } p+1
\\[0.5em]
	(\epsilon)	& \S^+\text{--$n$--recurrence} &	\S^+_p(n+1) = \S^+_p(n) + (n+1)^p
\\[0.5em]
	(\theta)	& \text{First $\S^+$ values} &	\S^+_p(1)=1; \S^+_p(0)=0; \S^+_p(-1)=-0^p
\\[0.5em]
	(\iota)	& \text{Definition}	&	\S^+_p(n) =\ldots+\B^+_p n, \text{ or }
								 \B^+_p := \tfr \S^+_p(n)/n\Big|_{n=0}
\\[0.5em]
	(\kappa)	& \B^+\text{--recurrence}	&	\binom{\B^+-1}{p} = \tfr {(-1)^p}/{(p+1)}
\\[0.5em]
	(\mu)	& \text{First numbers} &	\B^+_0=1; \B^+_1=\fr 1/2; \B^+_2=\fr 1/6; \B^+_3=0; \B^+_4=-\fr 1/{30}
\\[0.5em]
	(\nu)	& \text{Sums of powers}	&	\S^+_p(n) = \dfr {(\B^++n)^{p+1}-\B^+_{p+1}}/{p+1}
								:= \dfr {\sum_{k=0}^p\binom{p+1}{k}\B^+_kn^{p+1-k}}/{p+1}
									
\\[0.5em]
	(\xi)	& \text{Parity}			&	\B^+_{2k+1}=0, \text{ for } k>0		\phantom{\dfr 1/2}
\\
\hline
\end{array}
 $$

	Assuming the validity of the first table, the justification of this second table is immediate since, for $p>0$,
$$
	\S^+_p(n) = S_p(n) + n^p	\quad\implies\quad	\B^+_p=(-1)^pB_p.
 $$
The $\S^+$--$p$--recurrence ($\beta$) is the symbolic form of the one presented without proof by Bernoulli, and the $\B^+$--recurrence ($\kappa$) follows after dividing by $n$ and then setting $n=0$.

We note that from $(\gamma)$ and $(\theta)$, $2\S^+_1$ divides $\S^+_p$ if $p>0$.

\newpage

	Alternatively, Jacob Bernoulli could have used a simple induction to justify directly his general formula noted by $(\nu)$ in the preceeding table, should he have felt the need to provide a formal proof. Let us take his crucial rule for extending the list of numbers as a definition, and assume that for a fixed integer $p\ge0$, the numbers $\B^+_k, 0\le k\le p$ have been choosen so as to satisfy the formula $(\nu)$ for $n=1$ and the successive powers $0,1,2,\ldots,p$, so that
$$
\S^+_k(1) = \dfr {(\B^++1)^{k+1}-\B^+_{k+1}}/{k+1},	\qquad(0\le k\le p),
 $$
or \cite{Cesaro:1880a}, replacing $k+1$ by $k$, noting that $\S^+_k(1)=1$ and that empty sums are zero by convention,
$$
	k = (\B^++1)^{k}-\B^+_{k} = \sum_{j=0}^{k-1}\binom{k}{j}\B^+_j,	\qquad(0\le k\le p+1).
 $$
Now we can suppose by induction on $n$ that the formula $(\nu)$ has been verified for this fixed power $p\ge0$ and a certain $n\ge1$, and we consider the next case $n+1$.
\begin{align*}
(\B^+ + n + 1)^{p+1} - \B^+_{p+1} &= \left[(\B^+ + n)^{p+1} - \B^+_{p+1}\right]
							+ \left[(\B^+ + 1 + n)^{p+1} - (\B^+ + n)^{p+1}\right]
\\
	&= (p+1)\S^+_p(n) + \sum_{k=0}^{p+1} \binom{p+1}{k} \left[(\B^+ + 1)^k - \B^+_k \right] n^{p+1-k}
\\
	&= (p+1)\S^+_p(n) + \sum_{k=1}^{p+1} \binom{p+1}{k} \, k \, n^{p+1-k}
\\
	&= (p+1)\S^+_p(n) + (p+1)\sum_{j=0}^{p} \binom{p}{j} n^{p-j}
\\
	&= (p+1)\left[ \S^+_p(n) + (n+1)^{p}\right]
\\
	&= (p+1) \S^+_p(n+1).
\end{align*}
Thus the formula $(\nu)$ is also valid for $n+1$, which completes its proof by induction on the variable $n$ for a fixed integer power $p\ge0$, using the original definition of the Bernoulli numbers.

In summary, the Bernoulli general formula is valid for all $n\ge1$ if and only if it is valid for $n=1$, and we note that only the Newton binomial theorem was required here. A 2014 book by Arakawa, Ibukiyama, and Kaneko also adopts the $\B^+_1=\tfr 1/2$ convention but uses instead Taylor's theorem for polynomials for the proof of $(\nu)$, a calculus method dating back to an 1846 note by Arndt \cite{Arndt:1846a}.

Finally, we present the exponential generating function proof of formula $(J)$ for $S_p(n)$, based on the Euler definition of the Bernoulli numbers as Maclaurin coefficients of his very special function,
$$
	\fr x/{e^x-1} = \fr x/2\coth\fr x/2 - \fr x/2 = \sum_{k=0}^\infty B_k \fr x^k/{k!},
			\qquad(|x|<2\pi, B_1=-\fr 1/2).
 $$
This more advanced method elegantly deduces at once formula $(J)$ for all $p\ge0$ as follows.
\begin{align*}
	\sum_{p=0}^\infty S_p(n)\fr x^p/{p!} &= \sum_{p=0}^\infty \left\{\sum_{k=0}^{n-1} k^p\right\} \fr x^p/{p!}
	= \sum_{k=0}^{n-1} \sum_{p=0}^\infty  k^p \fr x^p/{p!}
	= \sum_{k=0}^{n-1} e^{kp} = \fr x/{e^x-1} \times \fr e^{nx}-1/x 
\\
	&= \left(\sum_{k=0}^\infty B_k\fr x^k/{k!}\right)\left(\sum_{j=0}^\infty \fr n^{j+1}/{j+1} \fr x^j/{j!}\right)
	 = \sum_{p=0}^\infty \left\{ \sum_{k=0}^p \binom{p}{k} B_k \fr n^{p+1-k}/{p+1-k} \right\} \fr x^p/{p!}
\\
	&= \sum_{p=0}^\infty \left\{ \sum_{k=0}^p \binom{p+1}{k} \fr B_k n^{p+1-k}/{p+1}\right\} \fr x^p/{p!}
	= \sum_{p=0}^\infty \left\{\fr (B+n)^{p+1}-B_{p+1}/{p+1}\right\} \fr x^p/{p!}.
\end{align*}
The proof of formula $(\nu)$ for $\S^+_p(n)$ would be similar, using the original Euler definition,
$$
	\fr x/{e^x-1} + x = \fr xe^x/{e^x-1} = \sum_{k=0}^\infty \B^+_k \fr x^k/{k!},
			\qquad(|x|<2\pi, \B^+_1=\fr 1/2).
 $$

\newpage


\section{Equivalent definitions}
Many different equations involve the Bernoulli numbers, in particular the even indexed $B_{2n}$, and some yield equivalent definitions  \cite{Gould:1972}. Currently, the first two in the following table are used most often (``E.G.F." means ``exponential generating function"). We use the Blissard representative notation for power series and current conventions with $B_1=-\tfr 1/2$, so that $\B^+_n=(-1)^nB_n$.

$$
\begin{array}{|l|l|l|} \hline
	\text{Method}	& \text{Equation}  & \text{References}
\\[0.5em] \hline
\text{Euler E.G.F.}	&	\e{Bx} := \sum_{n=0}^\infty B_n \tfr{x^n}/{n!}
												= \tfr x/{(\e{x}-1)}, \qquad(|x|<2\pi)			\phantom{\bigg[}
							& \text{\cite[p. 283]{Blissard:1861}}
\\[0.5em]
\text{De Moivre}		&	(B+1)^{p+1}-B_{p+1} := \sum_{k=0}^p\binom{p+1}{k}B_k = 0^p, \qquad(p\ge0)
							& \text{\cite[p. 280]{Blissard:1861}}
\\[0.5em]
\text{Blissard E.G.F.s} & 	 \cos 2Bx = x\cot x \text{ and } \sin 2Bx =-x, \qquad(|x|<\pi)
							& \text{\cite[p. 285]{Blissard:1861}}
\\[0.5em]
\text{Euler convolution} &	(2n+1)B_{2n} = - \sum_{k=1}^{n-1} \binom{2n}{2k}B_{2k}B_{2n-2k},\quad(n>1)
							& \text{\cite{Guinand:1979,Williams:1953}}
\\[0.5em]
\text{Plana integral} &	(-1)^{n-1}B_{2n} =
						4n\int_0^\infty \fr {t^{2n-1}}/{\e{2\pi t}-1}\,dt =
						\pi \int_0^\infty \fr {t^{2n}}/{\sinh^2(\pi t)}\,dt,\;(n>0)
							& \text{\cite{Plana:1820}}
\\[0.5em]
\text{Glaisher series} &	B_{4n+2} =	2(4n+2)\sum_{k=0}^\infty \fr {k^{4n+1}}/{e^{2\pi k}-1}	 + \tfr {0^n}/{4\pi}
							& \text{\cite{Glaisher:1889}}
\\[0.5em]
\text{Plana integral} &	(-1)^{n-1}(1-2^{1-2n})B_{2n} =
						\int_0^\infty \fr 4n{t^{2n-1}}/{\e{2\pi t}+1}\,dt =
						\int_0^\infty \fr \pi {t^{2n}}/{\cosh^2(\pi t)}\,dt
							& \text{\cite{Plana:1820}}
\\[0.5em]
\text{Glaisher series} &	(2^{4n+1}-1)B_{4n+2} = 2(4n+2)\sum_{k=1,\text{odd}}^\infty \fr {k^{4n+1}}/{e^{k\pi}+1}
							& \text{\cite{Glaisher:1889}}
\\[0.5em]
\text{Jensen integral} & (-1)^{n}B_{n} =	\fr \pi/2 \int_{-\infty}^\infty \fr {(\tfr 1/2+it)^n}/{\cosh^2(\pi t)}\,dt
							& \text{\cite{Touchard:1956}}
\\[0.5em]
\text{Riemann function} &	\zeta(2n) := \sum_{k=1}^\infty \tfr 1/{k^{2n}} = (-1)^{n-1}\tfr {(2\pi)^{2n}B_{2n}}/{(2(2n)!)}
							& \text{\cite{Euler:1755}}
\\[0.5em]
\text{Riemann function} &  (-1)^n B_n = -n\zeta(1-n)
							& \text{\cite{Knopp:1922,Sondow:1994}}
\\[0.5em]
\text{Blissard binomial}	&	\binom{B}{p} := B\ldots(B-p+1)/p! = \tfr {(-1)^p}/{(p+1)} 
							& \text{\cite[p. 284]{Blissard:1861}}
\\[0.5em]
\text{Blissard binomial} &  \binom{B+p-1}{p}  := B\ldots(B+p-1)/p! = - \tfr 1/{[p(p+1)]}, \quad(p\ge1)
							& \text{\cite[p. 284]{Blissard:1861}}
\\[0.5em]
\text{Blissard calculus} &  (-1)^nB_n = \sum_{k=0}^n \tfr {\sum_{j=0}^k (-1)^j\binom{k}{j} (j+1)^n}/{(k+1)}
							& \text{\cite[p. 67]{Blissard:1862}}		
\\[0.5em]
\text{Genocchi sum}	& 2^m(2^n-1)B_n = n \sum_{k=2}^m\binom{m}{k}\sum_{j=1}^{k-1}(-1)^{j-1}j^{n-1}, \;(m\ge n>1)
							& \text{\cite{Genocchi:1852}}
\\[0.5em]
\text{Matrix inverse}	&   \left[\, \tfr {B_{i-j}}/{(i-j)!} \,\right]_{n\times n,j\le i} =
								 \left[\, \tfr 1/{(i-j+1)!} \,\right]^{-1}_{n\times n,j\le i}
							& 
\\[0.5em]
\text{Determinant}	&  (-1)^n B_n = (n+1)!\,\left|\,\binom{i+1}{j-1}\,\right|_{n\times n,j\le i+1}
							& \text{\cite{Glaisher:1877,Hammond:1875}}
\\[0.6em]
\text{Determinant}	&  (-1)^n B_n = n! \,\left|\,\tfr 1/{(i-j+2)!}\,\right|_{n\times n,j\le i+1}
							& \text{\cite{Glaisher:1877,Sciacci:1865}}
\\[0.5em]
\hline
\end{array}
 $$
The ``Blissard calculus'' formula, not attributed in \cite[p. 82]{Saalschutz:1893}, gives a simple example of his symbolic manipulations. With $e^\theta=:1+x$ in his symbolic Euler E.G.F., he applies $n$ times the operator $\tfr -d/{d\theta}=-(1+x)\tfr d/{dx}$, and finally sets $x=\theta=0$. ``Required to express $B_n$ in direct terms''~:
\begin{align*}
e^{B\theta} = (1+x)^B &:= \fr \theta/{e^\theta-1} = \fr {\log(1+x)}/x  =  \sum_{k=0}^\infty \fr 1^0/{k+1} (-x)^k,		\qquad(|x|<1)
\\ \implies
-B(1+x)^B &:= \sum_{k=1}^\infty \fr k/{k+1}\,(1+x)(-x)^{k-1} = \sum_{k=0}^\infty \left(\fr 1^1/{k+1}-\fr {2^1-1^1}/{k+2}\right) (-x)^k
\\ \implies
(-B)^n(1+x)^B &:= \sum_{k=0}^\infty \left(\fr 1^n/{k+1}-\fr 2^n-1^n/{k+2}+\fr 3^n-2\cdot2^n+1^n/{k+3}- \ldots\right) (-x)^k.
\end{align*}
This uses the differences of $1^n$. The concise symbolic form $B_n=\{\fr \log(1+\Delta)/\Delta\}0^n$ is in \cite[p. 95]{Blissard:1867}.


\smallskip

\section{Conclusion}
We have shown that the original simple rule provided by Jacob Bernoulli to compute sequentially the Bernoulli numbers (via the De Moivre recurrence) is sufficient to prove many of the well known properties of this sequence of rational numbers, except their alternating signs.

\newpage


\end{document}